\documentclass[12pt]{article}
\usepackage{amsmath,amssymb,amscd,latexsym,mathrsfs}

\usepackage[all]{xy}
\usepackage[cp1251]{inputenc} 
\usepackage[english]{babel}

\newcommand{\colim}{{\rm colim}}
\newcommand{\cocon}{{\rm con\,}}
\newcommand{\hocolim}{{\rm hocolim}}
\newcommand{\coLim}{\underrightarrow{\lim}}

\newcommand{\Cat}{{\rm Cat}}
\newcommand{\Set}{{\rm Set}}
\newcommand{\sSet}{{\rm sSet}}
\newcommand{\psSet}{{\rm Set^{\Delta^{op}}_{\star}}}

\newcommand{\Grp}{{\rm Grp}}
\newcommand{\Ab}{{\rm Ab}}
\newcommand{\Ob}{{\rm Ob}}

\newcommand{\Imm}{{\rm Im\,}}
\newcommand{\Ker}{{\rm Ker\,}}
\newcommand{\Coker}{{\rm Coker\,}}
\newcommand{\diag}{{\rm diag\,}}

\newcommand{\NN}{{\,\mathbb N}}

\newcommand{\mC}{{\,\mathscr C}}
\newcommand{\mD}{{\,\mathscr D}}

\newcommand{\mG}{{\cal G}}

\newcommand{\ZZ}{{\,\mathbb Z}}

\newcommand{\mA}{{\mathscr A}}

\newcommand{\cA}{{\mathcal A}}

\newtheorem{theorem}{\bf Theorem}[section]
\newtheorem{lemma}[theorem]{\bf Lemma}
\newtheorem{proposition}[theorem]{\bf Proposition}
\newtheorem{corollary}[theorem]{\bf Corollary}
\newtheorem{definition}{\sc Definition}[section]
\newtheorem{example}[definition]{\sc Example}
\newtheorem{remark}[definition]{\sc Remark}

\def\leq{\leqslant}
\def\geq{\geqslant}

\def\MYstar{\mathop{\star}\limits}

\begin{document}

\begin{center}
{\large \bf Non-Abelian homology 
and homotopy colimit\\
 of classifying spaces for a diagram of groups\\
} 
\medskip
Ahmet A. Husainov\\
\end{center}

\begin{abstract}
This paper considers non-Abelian homology groups of a group diagram introduced as homotopy groups of a simplicial change.
We prove a theorem stating that the non-Abelian homology groups of a group diagram are isomorphic to the homotopy groups of the homotopy colimit of a classifying space diagram, with the dimension shifted by 1.
Bousfield and Kan proved an isomorphism between the homotopy groups of an Abelian simplicial group and the homology groups of this simplicial group.
We generalize this to non-Abelian simplicial groups.
We also develop a method for finding a non-zero homotopy group of smallest dimension for the homotopy colimit of classifying spaces.
For a group diagram over a free category with a zero colimit, we obtain a criterion for the isomorphism of the first non-Abelian and Abelian homology groups.
\end{abstract}

2000 Mathematics Subject Classification 55U10, 18C15, 18G10, 55P10, 55R35


\tableofcontents

\section{Introduction}

We introduce non-Abelian homology $\colim^{\mC}_n \mG$ for an arbitrary group diagram $\mG$ over the small category $\mC$ (Definition \ref{defhomol}) using a simplicial replacement of the group diagram in the sense of Bousfield and Kan \cite[XII.5.4]{bou1972}.
We show that these homology groups can be interpreted as homotopy groups of the homotopy colimit for the diagram of classifying spaces $B\mG$ 
(Theorem \ref{main2}).

In homotopy theory, many good methods and models have appeared for the study and calculation of homotopy groups \cite{bro2011}.
However, it remains necessary to calculate the homology,
and this requires the development of homological methods.
This paper is devoted to the theory of non-Abelian homology of group diagrams,
 which develops the classical theory of homology of small categories.
The theory of non-Abelian homology of diagrams is aimed at solving problems of homotopy theory.

The homology of small categories with coefficients in diagrams of abelian groups, shortly abelian homology of diagrams, can be defined in a similar way, using a simplicial replacement.
However, the meanings of non-Abelian and Abelian diagram homologies are very different. For example, for the category $V$ consisting of three objects and two morphisms
$$
b\leftarrow a\rightarrow c,
$$
  except for the identical ones, the Abelian homology of any diagram over $V$ is equal to $0$
   in dimensions $n\geq 2$.
Let's compare them with non-Abelian homology.
To this aim, we use the following formula, proved in our paper (Theorem \ref{main2}):
\begin{equation}\label{mainform}
\colim^{\mC}_n \mG\cong \pi_{n+1}(\hocolim^{\mC}B\mG),
\end{equation}
where $\hocolim^{\mC}B\mG$ is the homotopy colimit of the diagram in the category of pointed simplicial sets, which consists of the classifying spaces (nerves) of the groups of the diagram 
$\mG$.
It follows from this formula that the $n$th non-Abelian homology of the group diagram $0\leftarrow \ZZ\to 0$ is isomorphic to the homotopy groups of the two-dimensional sphere 
$\pi_{n+1}(S^2)$.
It was established in \cite{iva2016} that the $n$th homotopy groups of the two-dimensional sphere are not equal to $0$ for all $n\geq 2$. This implies that the $n$th non-Abelian homology of the diagram
$0\leftarrow \ZZ\to 0$ are not equal to its Abelian homology for every $n\geq 1$.

The formula (\ref{mainform}) also shows that the non-Abelian homology of the diagram $0\leftarrow G \to 0$ is isomorphic to the homotopy groups of the suspension over
 the classifying space of the group $G$, which in some cases were calculated in the papers \cite{bro1987} and \cite{mik2010}.
In \cite{rom2013}, the software for calculation of homotopy groups of such spaces is developed. This points to the possibility of computing non-Abelian homology of group diagrams.

 The assertion that we obtained from the formula (\ref{mainform}) that for simplicial groups considered as group diagrams (over $\Delta^{op}$), non-Abelian homology is isomorphic to homotopy groups also inspires hope for perspective (Corollary \ref{cor53}).
 
From the formula (\ref{mainform}) we also obtain a criterion for the simply connectedness 
of the homotopy colimit for the diagram of classifying spaces of groups (Corollary \ref{connect}).
By the formula (\ref{mainform}), the problem of calculating non-Abelian homology of a group diagram
\begin{equation}\label{prob22}
\mG(b) \leftarrow \mG(a) \rightarrow \mG(c)
\end{equation}
is closely related to some problems of describing the homotopy groups of a homotopy pushout, for example \cite[Problem 2.2]{bro1990}.
Note also that the formulas for the homotopy groups of the homotopy colimit of the diagram of classifying spaces of quotient groups $G\diagup N_{i_1}\cdots N_{i_k}$ over the inclusion-partially ordered set of proper subsets of the set $\{1, 2, \ldots, n\}$ proved by Ellis and Mikhailov \cite{ell2010} can be considered as formulas for non-Abelian homology in $n$ and $n-1$ dimensions.

In our paper, for a group diagram $\mG$ over a free category $W\Gamma$ with 
$\colim^{W\Gamma}\mG=0$, we obtain an isomorphism criterion for the first non-Abelian homology group of this diagram and the first Abelian homology group of its abelianization (Proposition \ref{compar}).

We have introduced the connectivity $\cocon\mG$ of the group diagram $\mG$ as the smallest $n\geq 0$ for which the group $\colim^{\mC}_n\mG$ is non-trivial. Example \ref{expar2} demonstrates a method for finding the connectivity $n= \cocon {\mG}$ and calculating 
the group $\colim^{\mC}_n {\mG}$. The number $\cocon\mG$ is equal to the largest $n$ for which the space $\hocolim^{\mC}B\mG$ is $n$-connected.

 \section{Preliminaries}
 
  We work with simplicial sets. We adhere to the notation from Gabriel and Zisman \cite{gab1967} and MacLane \cite{mac2004}.
Homotopy colimits will be considered in the category of pointed simplicial sets, the main source is the monograph by Bousfield and Kan \cite{bou1972}.
Pointed bisimplicial sets will also be used, from the book by Goerss and Jardine \cite{goe2009} and some information about simplicial groups and their homotopy groups from Wu \cite{wu2010}.

\subsection{Notation}

Let us write out the notation and some definitions. The rest will be given in the course of the presentation.

\begin{itemize}
\item
$\mA^{\mC}$ is the category of functors from the small category $\mC$ to a 
category $\mA$.
 Functors $F$ from a small category $\mC$ to an arbitrary $\mA$ are called object diagrams of the category $\mA$ over $\mC$ and can be denoted as the family $\{F(c)\}_{c\in \mC}$.
 \item $\Delta{A}=\Delta_{\mC}A$ is a functor $\mC\to \mA$ that takes constant values 
 $A\in \Ob\mA$ on objects and $1_A$ on morphisms.
\item $\NN$ - set of non-negative integers.
\item $\Cat$ - category of small categories and functors.
\item
$\Set$ - category of sets and mappings.
\item
${\Grp}$ is the category of groups and homomorphisms.
\item $0$ - a group consisting of one element.
\item $\Delta$ is the category of finite linearly ordered sets $[n]= \{0, 1, \cdots, n\}$, $n\geq 0$, and non-decreasing mappings. The category $\Delta$ is generated by morphisms of the following form:
\begin{enumerate}
\item $\partial^i_n: [n-1]\to [n]$ (for $0\leq i\leq n$) - increasing mapping whose image does not contain $i$,
\item
$\sigma^i_n: [n+1]\to [n]$ (for $0\leq i\leq n$) is a non-decreasing surjection that takes the value $i$ twice,
\end{enumerate}
\item
$\colim^{\mC}_n: \Grp^{\mC}\to \Grp$ ($\forall n\geq 0$) are functors assigning to each group diagram the values of the non-Abelian homology of this diagram (the definition will be given in our work),
\item
$\coLim^{\mC}_n: \Ab^{\mC}\to \Ab$ ($\forall n\geq 0$) are Abelian homology functors defined in \cite[Application II]{gab1967} as left satellites of the colimit functor, and in \cite{X2022} as homology of small categories with coefficients in object diagrams of the abelian category.
\end{itemize}

 Let $\Phi: \mC\to \mD$ be a functor between small categories.
For an arbitrary category $\mA$, the inverse image functor $\Phi^*: \mA^{\mD}\to \mA^{\mC}$ is defined, which assigns to each diagram $F\in \mA^{\mD}$ composition $F\Phi= F\circ \Phi\in \mA^{\mC}$, and to each natural transformation $\eta: F\to F'$ the natural transformation $\eta\Phi: F\Phi \to F'\Phi$ defined by the formula $(\eta\Phi)_{c}= \eta_{\Phi(c)}$, for all $c\in \Ob\mC$.
If $\mA$ is a cocomplete category, then the functor $\Phi^*$ has a left adjoint functor $Lan^{\Phi}: \mA^{\mC}\to \mA^{\mD}$, which is called the left Kan extensions.

  A simplicial set is a functor $X: \Delta^{op}\to \Set$. The values of this functor on morphisms generating the category $\Delta$ are denoted by $d^n_i= X(\partial^i_n)$, $s^n_i= X(\sigma^i_n)$, for all $0\leq i\leq n$.
A simplicial mapping $X\to Y$ between simplicial sets is a natural transformation.

A simplicial set $X$ is pointed if a vertex $x\in X_0$ is distinguished in it. 
We call the distinguished vertex $x$ to be base point. 
A morphism of simplicial sets is said to be pointed if it maps a base point to a base one.

The category of simplicial sets is denoted by $\Set^{\Delta^{op}}$ or $\sSet$. The category of pointed simplicial sets and pointed morphisms is denoted by $\sSet_*$.

The classifying space (nerve) of the small category $\mC$ is denoted by $B\mC$.

 Consider the category of simplices $(\Delta\downarrow X)^{op}$ for the simplicial set $X$ 
 dual to the comma category (in the sense of \cite[Chapter 2, \S2.6]{mac2004}) over $X$ with respect to the Yoneda embedding functor $\Delta(-): \Delta\to \Set^{\Delta^{op}}$.
Below we will see (Proposition \ref{systosym}) that for an arbitrary category with coproducts $\mA$ to every diagram $F: (\Delta\downarrow X)^{op}\to \mA$ there corresponds a simplicial object in $\mA $ consisting of coproducts $\coprod_{x\in X_n}F(x)$.
This object is used to construct a simplicial replacement of the diagram $F: \mC\to \mA$.

We begin our presentation by examining the following question.
There are two methods for constructing a simplicial replacement.
  Bousfield and Kan \cite{bou1972} define it as a simplicial object
  $\coprod\limits_{c_0\leftarrow c_1\leftarrow \ldots \leftarrow c_n} F(c_n)$.
  Quillen \cite{qui1973} uses the notation $\coprod\limits_{c_0\to c_1\to \ldots \to c_n} F(c_0)$.
  We will show that these two methods lead to isomorphic simplicial objects.

\subsection{Diagrams on the category of simplices}
 
Consider the category of simplices of a simplicial set.
Each diagram over this category in the category with coproducts $\mA$ corresponds to a simplicial object in $\mA$. We describe its face and degeneration operators.
 
\begin{definition}
  Let $X: \Delta^{op}\to \Set$ be a simplicial set.
Its category of simplices consists of the set of objects $\coprod\limits_{n\geq 0}X_n$. Morphisms of the category of simplices
$(m, x)\xrightarrow{\alpha} (n,y)$ are given by
nondecreasing mappings $\alpha: [n]\to [m]$ such that $X(\alpha)(x)= y$.
The composition of morphisms is
$$
((n, y)\xrightarrow{\beta} (p,z))\circ ((m, x)\xrightarrow{\alpha} (n,y))
= ((m, x)\xrightarrow{\alpha\beta} (p,z)).
$$
\end{definition}
 
The category of simplices is isomorphic to $(\Delta\downarrow X)^{op}$, and we will denote it by the same name.

We need the following statement about diagrams $(\Delta\downarrow X)^{op}\to \mA$ for the category $\mA$ with coproducts.

Let $Q^{op}_X: (\Delta\downarrow X)^{op}\to \Delta^{op}$ denote the forgetful functor defined as $(m, x)\mapsto [m]$ on objects and $(m, x)\xrightarrow{\alpha} (n, y) \mapsto \alpha: [n]\to [m]$ - on morphisms.
  Consider the left Kan extension $Lan^{Q^{op}_X}: \mA^{(\Delta\downarrow X)^{op}}\to \mA^{\Delta^{op}}$ along the functor $Q ^{op}_X$. We follow the definition of the left Kan extension given in \cite[\S 10.3]{mac2004}, with the values of the left Kan extension being computed pointwise.
The following assertion is proved in exactly the same way as the first part of 
\cite[Proposition 3.2]{X2022} for the case $\mD= \Delta$.

\begin{proposition}\label{systosym}
For any diagram $F: (\Delta\downarrow X)^{op}\to \mA$ into the category $\mA$ with coproducts, the value of $Lan^{Q^{op}_X}F$ is equal to the simplicial object consisting from objects $Lan^{Q^{op}_X}F [n]= \coprod\limits_{x\in X_n}F(x)$.
The boundary operators $d^n_i$ ($n\geq 1$) and degeneracy $s^n_i$ ($n\geq 0$) of this simplicial object are defined as morphisms making the following diagrams commutative for all $x\in X_n$ and $0 \leq i\leq n$:
$$
\xymatrix{
\coprod_{x\in X_n}F(x) \ar[r]^{d^n_i} & \coprod_{x\in X_{n-1}}F(x)\\
F(x)\ar[u]^{in_x} \ar[r]_{F(x\xrightarrow{\partial^i_n}X(\partial^i_n)x)}
  & F(X(\partial^i_n)x)\ar[u]_{in_{X(\partial^i_n)x}}
}\quad
\xymatrix{
\coprod_{x\in X_n}F(x) \ar[r]^{s^n_i} & \coprod_{x\in X_{n+1}}F(x)\\
F(x)\ar[u]^{in_x} \ar[r]_{F(x\xrightarrow{\sigma^i_n}X(\sigma^i_n)x)}
  & F(X(\sigma^i_n)x)\ar[u]_{in_{X(\sigma^i_n)x}}
}
$$
\end{proposition}

 \begin{example}
If $\mA= {\Grp}$, then for an arbitrary group diagram $F$ over the category of elements $(\Delta\downarrow X)^{op}$ we obtain a simplicial group consisting of free products of the groups $\MYstar_{x\in X_n}F(x)$.
The boundary and degeneration operators can be written explicitly
\begin{enumerate}
\item $d^n_i (x\in X_n, f\in F(x))= (X(\partial^i_n)x, \label{opbd1}
F(x\xrightarrow{\partial^i_n}X(\partial^i_n)x) (f))$, \label{opbd2}
\item
$s^n_i (x\in X_n, f\in F(x))= (X(\sigma^i_n)x,
F(x\xrightarrow{\sigma^i_n}X(\sigma^i_n)x) (f))$.
\end{enumerate}
  \end{example}

 \subsection{Simplicial replacement}
 
  The classifying space of a small category $\mC$ is its nerve, the simplicial set $B\mC =\Cat(-, \mC)|_{\Delta}: \Delta^{op}\to \Set$.
   Its $n$-dimensional simplices $x\in B_n\mC$ can be viewed as functors $x: [n]\to \mC$ or sequences of morphisms $c_0\xrightarrow{\alpha_1} \ldots \xrightarrow{\alpha_n} c_n $. The boundary operators $B\mC(\partial^i_n): B_n\mC\to B_{n-1}\mC$ map $x\in \Cat([n], \mC) \mapsto x\circ \partial ^i_n$ or, equivalently, they map to the tuple $(c_0\xrightarrow{\alpha_1} \ldots \xrightarrow{\alpha_n} c_n)$, $c_i=x(i)$, one of the sequences
   \begin{itemize}
   \item $(c_0\xrightarrow{\alpha_1} \ldots
   \xrightarrow{\alpha_{n-1}} c_{n-1})$ if $i=n$,
   \item $(c_0\xrightarrow{\alpha_1} \ldots
   \to c_{i-1} \xrightarrow{\alpha_{i+1}\alpha_{i}} c_{i+1} \to \ldots
   \xrightarrow{\alpha_n} c_n)$ if $0< i< n$,
   \item $(c_1\xrightarrow{\alpha_2} \ldots
   \xrightarrow{\alpha_n} c_n)$ if $i=0$.
   \end{itemize}
  The degeneration operators $B\mC(\sigma^i_n): B_n\mC\to B_{n+1}\mC$ can be calculated (for all $0\leq i\leq n$) by the formula $x\mapsto x\circ \sigma^i_n$ or by
  $$
  (c_0\xrightarrow{\alpha_1} \ldots \xrightarrow{\alpha_n} c_n) \mapsto (c_0\xrightarrow{\alpha_1} \ldots \xrightarrow{\alpha_{i}} c_i \xrightarrow{1_{c_i}} c_i\xrightarrow{\alpha_{i+1}} \ldots \xrightarrow{\alpha_n} c_n).
  $$
  This implies that the categories $\Delta\downarrow \mC$ and $\Delta\downarrow B\mC$ are naturally isomorphic (with respect to $\mC\in \Cat$), and we can consider them as equal categories.

Consider the functors
  $$
  \mC\xleftarrow{\partial_0} (\Delta\downarrow B\mC)^{op}
  \xrightarrow{Q^{op}_{B\mC}} \Delta^{op},
  $$
  the first of which maps each $x:[m]\to \mC$ corresponding to some tuple $(c_0\xrightarrow{\alpha_1} \ldots \xrightarrow{\alpha_m} c_m)$ to its object $x(0)= c_0$ and to each morphism $x\to y$ of the category $(\Delta\downarrow \mC)$:
$$
\xymatrix{
[n] \ar[rr]^{y} && \mC\\
& [m]\ar[lu]^g \ar[ru]_x
}
$$
morphism $y(0) \xrightarrow{y(0\leq g(0))} y(g(0))=x(0)$ of the category $\mC$. This functor $\partial_0$ is contravariant.
 
  The functor $Q^{op}_{B\mC}$ is equal to the functor $Q^{op}_X$ defined above for $X=B\mC$.

\begin{definition}\label{defrepl}
  The simplicial replacement functor $\coprod: \mA^{\mC}\to \mA^{\Delta^{op}}$ is equal to the composition of the left Kan extension $Lan^{Q^{op}_{B\mC}}: \mA^{(\Delta\downarrow B\mC)^{op}} \to \mA^{\Delta^{op}}$ and the functor $\partial^*_0: \mA^{\mC} \to \mA^{(\Delta\downarrow B\mC)^{op}}$ defined by $(F\in \mA^{\mC})\mapsto 
  F\circ \partial_0$.
   In particular, the simplicial change of any diagram $F$ over $\mC$ is equal to $Lan^{Q^{op}_{B\mC}}(F\partial_0)$.
\end{definition}

Substituting the composition $F\partial_0$ into the sentence \ref{systosym} instead of $F$
 and $X= B\mC$, we arrive at the following assertion.
  
\begin{proposition}
The simplicial replacement of the diagram $F: \mC\to \mA$ consists of the objects $\coprod_n F= \coprod\limits_{c_0\to c_1 \to \ldots \to c_n}F(c_0)$.
Boundary operators are defined as morphisms $d^n_i: \coprod_n F\to \coprod_{n-1}F$, $0\leq i\leq n$ making commutative the following diagram on the left for every $i$ in the range $0<i \leq n$, and for $i=0$ - on the right:
$$
\xymatrix{
\coprod_n F \ar[r]^{d^n_i} & \coprod_{n-1}F\\
F(c_0)\ar[u]|-{in_{c_0\to \ldots \to c_n}} \ar[r]_{=}
  & F(c_0)\ar[u]_{in_{B\mC(\partial^i_n){(c_0\to \ldots \to c_n)}}}
} \qquad
\xymatrix{
\coprod_n F \ar[r]^{d^n_0} & \coprod_{n-1}F\\
F(c_0)\ar[u]|-{in_{c_0\to \ldots \to c_n}} \ar[r]_{\alpha_1}
  & F(c_1)\ar[u]|-{in_{(c_1\to \ldots \to c_n)}}
}
$$
The degeneration operators $s^n_i: \coprod_n F \to \coprod_{n+1}F$ are defined as morphisms making the following diagrams commutative for all $0\leq i\leq n$:
$$
\xymatrix{
\coprod_n F \ar[r]^{s^n_i} & \coprod_{n+1}F\\
F(c_0)\ar[u]|-{in_{c_0\to \ldots \to c_n}} \ar[r]_{=}
  & F(c_0)\ar[u]_{in_{B\mC(\sigma^i_n){(c_0\to \ldots \to c_n)}}}
}
$$
\end{proposition}

Bousfield and Kan \cite{bou1972} consider the classifying space of the dual category $X= B\mC^{op}$ and obtain the simplicial replacement $\coprod' F = Lan^{Q^{op}_{B\mC^ {op}}}F\partial$, where $\partial(c_0\leftarrow c_1 \leftarrow \ldots \leftarrow c_n)=c_n$.
Since the simplicial sets $B\mC$ and $B\mC^{op}$ are isomorphic, the resulting simplicial replacements are isomorphic.

\subsection{Homotopy colimit of pointed simplicial
sets}

A simplicial set $X$ is said to be pointed if a point $\star_0\in X_0$ is distinguished in it.
The distinguished point is called base point.
For an arbitrary simplicial set $X$, we can select any point $\star_0$ from the set $X_0$.
Together with the point $\star_0\in X_0$, for each $n \geq 0$ there will be a single simplex 
$\star_n= X([n]\to [0]) \star_0)\in X_n$.
Since every diagram in the category $\Delta$ consisting of morphisms
$$
\xymatrix{
[0] && [m ]\ar[ll] \ar[dl]^f\\
&[n]\ar[ul]
}
$$
is commutative, then faces and degeneracies leave simplices in this sequence, and we can consider a pointed simplicial set as a functor from $\Delta^{op}$ to the category of pointed sets whose values are equal to $(X_n, \star_n)$.
The category of pointed simplicial sets is isomorphic to the category of functors $\Delta^{op}\to \Set_{\star}$, but homotopy groups are considered with respect to the base point 
$\star_0$.

Let $\sSet= \Set^{\Delta^{op}}$ be the category of simplicial sets, and $\sSet_{\star}$ be the category of pointed simplicial sets and simplicial mappings that preserve the base point.

The homotopy colimit of a diagram of pointed simplicial sets is well described by Bousfield and Kan \cite[Chapter XII]{bou1972}.
In \cite[Lemma 5.2]{bou1972} it is proved that the homotopy colimit functor $\hocolim^{\mC}:(\sSet_{\star})^{\mC} \to \sSet_{\star}$ can be defined as functor composition
$$
(\sSet_{\star})^{\mC} \xrightarrow{\coprod} (\sSet_{\star})^{\Delta^{op}}=
\Set^{\Delta^{op}\times \Delta^{op}}_{\star}
\xrightarrow{\diag} \psSet= \sSet_{\star},
$$
where $F\mapsto \coprod F= Lan^{Q^{op}_{\mC}} F\partial_0$ is the simplicial replacement functor for the case $\mA= \sSet_{\star}$ and $\diag$ - the diagonal functor \cite[Chapter XII]{bou1972}, corresponding to each bisimplicial set consisting of the sets 
$$
\{X([n],[p])\}_{([n],[p])\in \Delta^{op}\times \Delta^{op}}
$$ 
the simplicial set $\{X([n],[n])\}_{([n])\in \Delta^{op}}$.

Let $F: \mC\to \sSet_{\star}$ be an arbitrary diagram of pointed simplicial sets. Its simplicial replacement is equal to a bisimplicial set whose vertical (in the sense of 
\cite[Page 209]{goe2009}) simplicial sets are pointed and equal for each $n\geq 0$ to the coproduct (and hence to the wedge) $\bigvee\limits_{c_0 \to c_1 \to \ldots \to c_n} F(c_0)$ over all sequences of morphisms of length $n$.

Therefore, the homotopy colimit, as the diagonal of the resulting simplicial set, is equal to the pointed simplicial set
$$
(\hocolim^{\mC}F)_n = \bigvee\limits_{c_0\to c_1 \to \ldots \to c_n} F(c_0)_n, \text{for all $n\geq 0$}.
$$
  
\begin{example}\label{forconst}
Let $\Delta_{\mC}Y: \mC\to \sSet_{\star}$ be a diagram taking constant values $Y$. In order to construct $\hocolim^{\mC} \Delta_{\mC}Y$, we first construct a homotopy colimit of the diagram in which the base points are treated as ordinary.
  Let $U: \sSet_{\star}\to \sSet$ be a functor forgetting base points.
  Get
  $$
  (\diag\coprod \Delta_{\mC}Y)_n= \coprod_{c_0\to c_1\to \ldots \to c_n}Y_n =
   B_n\mC\times Y_n.
  $$
  Hence the simplicial set $\hocolim^{\mC}U(\Delta_{\mC}Y)$ is equal to $B\mC\times Y$.

Returning to $\hocolim^{\mC}\Delta_{\mC}Y$, we find that it consists of pointed sets 
$$
\bigvee\limits_{c_0\to c_1 \to \ldots \to c_n} Y_n
$$ 
equal to the sets $B_n\mC\times Y_n$ with identified simplices from $B_n\mC\times \star_n$. Therefore, the homotopy colimit of a diagram on $\mC$ taking constant values equal to the pointed simplicial set $Y$ is equal to the pointed simplicial factor set
$$
(B\mC\times Y)\diagup (B\mC\times \star).
$$
\end{example}

\section{Classifying space of a simplicial group}

This section is devoted to the homotopy groups of the simplicial group.
To every simplicial group there corresponds a diagram of nerves over $\Delta^{op}$. An auxiliary lemma is considered that the $n$th homotopy groups of the homotopy colimit of the nerve diagram are isomorphic to the $(n-1)$th homotopy groups of the simplicial group for all $n\geq 1$..

We adhere to the definition of the homotopy group $\pi_n(X)$ of a pointed simplicial set given in \cite[Chapter VI, \S 3]{gab1967}.
In particular, when defining the homotopy groups of a simplicial group its group structure can be ignored. But it is a fibrant simplicial set, and surjective morphisms of simplicial groups are (Kan)
fibrations.

\subsection{Classifying space of a discrete group}

If we consider a group as a small category, then we will define for it
its nerve --
a simplicial set consisting of the sets
$G^{[n]}= \Cat([n], G)$ and the mappings $G^{[n]}\to G^{[m]}$ corresponding
to the mappings $[m]\to [n]$. We denote this nerve by $BG$.
Let $EG$ be the simplicial set consisting of the sets
$E_n G= \{(g_0, \ldots, g_n) | (\forall i\in [n]) g_i\in G\}$.
It is the nerve of the category whose objects are $g\in G$,
and for each pair of objects $(g,g')$ there is a unique morphism
$g\to g'$.
We consider the homotopy groups of a simplicial set
as the homotopy groups of its geometric realization.
The functor from the category of zero morphisms to $G$ yields a simplicial
map of nerves $EG\to BG$, which is a bundle with discrete
fiber $G$. This
bundle
gives an exact sequence leading to
the formulas
$$
\pi_n(BG)\cong
\begin{cases}
1, & \text{ if $n=0$;}\\
G, & \text{ if $n=1$;}\\
0, & \text{ for all $n>1$.}
\end{cases}
$$
(Here $1$ denotes the set consisting of
a single element.)
This fact allows us to call the simplicial set $BG$
the {\it classifying space} of $G$.

\subsection{Moore complex and its homology groups}

Let $G: \Delta^{op}\to {\Grp} $ be a simplicial group.

We describe the Moore complex \cite[Ch.II, \S3.7]{qui1967} for calculating the homotopy 
groups $\pi_n(G)$.
For this purpose, consider the complex consisting of the groups $N_n G= \bigcap\limits_{i>0} \Ker (d^n_i: G_n\to G_{n-1})$ for $n>0$, and $ N_0 G=G_0$.

Its differentials $d^n: N_n G \to N_{n-1}G$ are defined by the restriction of homomorphisms $d^n_0$ for $n>0$, and $d^0=0$.

For example, for a constant simplicial group $G= \Delta G_0$, this gives isomorphisms 
$\pi_0(G)\cong G_0$ and $\pi_n(G)= 0$ for all $n>0$.

Note that the groups $\pi_n(G)$ (see \cite[Exercise 8.3.2]{wei1994}
or \cite[Chapter 3, Theorem 1.5]{wu2010}) are Abelian for $n\geq 1$, 
 while for $n=0$ they can be non-commutative.

The homotopy groups $\pi_n(G)$ are isomorphic to the homology groups of this complex $(NG,d_0)$
 (see \cite[Chapter 3, Theorem 1.5]{wu2010}):
$$
H_n(NG)= \Ker (d^n: N_n G\to N_{n-1}G )\diagup Im (d^{n+1}: N_{n+1}G\to N_n G).
$$

According to \cite[Chapter 3, Theorem 1.5]{wu2010}, Moore proved that for every simplicial
group $G$, the formulas
$$
H_n(NG, d_0) \cong \pi_n(G) \cong \pi_n(|G|)
$$
hold.
Since $C_*(\mC, \mG)$ is a simplicial group, this yields
\begin{corollary}\label{twoiso} For every diagram of groups $\mG: \mC\to \Grp$, there are isomorphisms
$$
\pi_n(C_*(\mC, \mG)) \cong \pi_n(|C_*(\mC, \mG)|) \cong H_n(NC_*(\mC, \mG))
$$
\end{corollary}

\subsection{Bisimplicial set and diagonal}

Consider an arbitrary bisimplicial set $X(-,=)$.
Let $\diag X(-,=)$ denote the diagonal of a double simplicial set.
It is a simplicial set defined as $\diag X([n])= X([n], [n])$.
A bisimplicial set $X(-,=)$ is called pointed if for every $m\geq 0$ the simplicial set $X(m,=)$ is pointed, i.e. has a marked vertex $x_m \in X(m,0)$, and for each morphism $\alpha \in \Delta([m], [n])$ the equality $X(\alpha)(x_n)= x_m$. Any bisimplicial set becomes pointed if we distinguish an element $x_0\in X(0,0)$ in it.

Let $G:\Delta^{op}\to {\Grp}$ be a simplicial group.
Consider the bisimplicial set $G(-)^{[=]}$ formed from the nerves of the groups and consisting of the sets $G(m)^{[n]}= B_n(G(m))$. We will consider it as a pointed simplicial set whose base point is the only element of the set $B_0(G(0))$

\begin{lemma}\label{spquill}
Let $G:\Delta^{op}\to {\Grp}$ be a simplicial group.
Then for all $n\geq 1$ there are natural isomorphisms
$$
\pi_n \diag(BG)
\cong
\pi_n(\diag G(-)^{[=]})
\cong\pi_{n-1}G.
$$
\end{lemma}
{\sc Proof.}
For each $p\geq 0$ the classifying space of the group $G(p)$ is connected. We choose a base point equal to the unique $0$-dimensional simplex of the space $BG(p)$ for each $p\geq 0$.
One can use the Bousfield-Friedlander spectral sequence \cite{bou1978} or the Zisman spectral sequence (see \cite[Appendix A]{bro1987}) for the bisimplicial set
$G(-)^{[=]}$:
$$
  E^2_{p,q} = \pi^h_p \pi^v_q G(-)^{[=]} \Rightarrow \pi_{p+q} (\diag G(-)^{[=]}) .
$$
For every $p\geq 0$ the groups $\pi_q G(p)^{[=]}$ are isomorphic to $G(p)$ for $q=1$, and trivial for $q\not=1$. This leads to degeneration of the spectral sequence and isomorphisms $E^2_{p,1}\cong \pi_{p+1} \diag G(-)^{[=]}$.
Hence $\pi_p(G)\cong \pi_{p+1}(\diag G(-)^{[=]})$ for all $p\geq 0$.
\hfill$\Box$
\begin{remark}
There are other methods for proving Lemma \ref{spquill}, for example, using the exact sequence of homotopy groups for the fibration $G \to \diag EG \to \diag BG$ leading to the natural isomorphism $\pi_n \diag(BG)\cong \pi_ {n-1}G$
\cite[Page 244]{goe2009}.
\end{remark}

\section{Non-Abelian homology with coefficients in a group diagram}

In this section, non-Abelian homology of a small category with coefficients in the group diagram is introduced. They are constructed as homotopy groups of a simplicial replacement of the group diagram. This replacement is a simplicial group, and hence the $n$th homology 
 for all $n\geq 0$ are groups, commutative for $n\geq 1$.

\subsection{Simplicial replacement of group diagram}

For an arbitrary small category $\mC$, the following functors are defined
$$
\mC \stackrel{\partial_0}\leftarrow (\Delta\downarrow \mC)^{op}\stackrel{Q_{\mC}^{op}}\to
\Delta^{op}.
$$

Let $\mA$ be a category with coproducts.
For each diagram $F: \mC\to \mA$ its simplicial replacement $\coprod F \in \mA^{\Delta^{op}}$ is defined as $\coprod F= Lan^{Q_{\mC}^ {op}}(F\circ\partial_0)$.
The object of $n$-dimensional simplices of the simplicial object $\coprod F$ will be equal to $(\coprod F)_n = \coprod\limits_{c_0\to \cdots \to c_n}F(c_0)$.
In particular, for the category of groups ${\Grp}$, it will be equal to the free product $\star_{c_0\to \cdots \to c_n} F(c_0)$, for the category of abelian groups $\Ab$ it will be equal to the direct sum $\oplus_ {c_0\to \cdots \to c_n} F(c_0)$, and in the case of the category of pointed simplicial sets $\sSet_{\star}$ it is equal to the wedge $\bigvee\limits_{c_0\to \cdots \to c_n} 
F (c_0)$.

\begin{proposition}\label{diagtosimp}
Let $\mG: \mC\to \Grp$ be a group diagram.
The simplicial group $\coprod \mG$ consists of the groups
$$
C_n(\mC, \mG)=\underset{c_0\to \cdots \to c_n}{\star}\mG(c_0).
$$
Its boundary operators $d^i_n: C_n(\mC,\mG)\to C_{n-1}(\mC,\mG)$ and degeneration operators $s^i_n: C_n(\mC,\mG)\to C_{n+1}(\mC,\mG)$, for $0\leq i\leq n$, act on the elements of the factors of free products by the formulas
\begin{multline}\label{fsimprep1}
d^i_n(c_0\xrightarrow{\alpha_1}c_1 \to
\cdots \to c_{n-1} \xrightarrow{\alpha_n} c_n, x)=\\
=
\begin{cases}
(c_1\xrightarrow{\alpha_2} c_2\to \cdots \to c_{n-1}
\xrightarrow{\alpha_n}c_n, \mG(\alpha_1)(x)),
& \text{ if $i=0$},\\
(c_0\xrightarrow{\alpha_1}c_1\to \cdots \to c_{i-1} \xrightarrow{\alpha_{i+1}\alpha_i}
c_{i+1}\to \cdots \xrightarrow{\alpha_n} c_n, x), & \text{ if $i>0$},
\end{cases}
\end{multline}
\begin{multline}\label{fsimprep2}
s^i_n(x, c_0\xrightarrow{\alpha_1}c_1 \to \cdots \to c_{n-1} \xrightarrow{\alpha_n} c_n)=\\
(c_0\xrightarrow{\alpha_1}c_1 \to \cdots \to c_i \xrightarrow{id} c_i
  \to \cdots \to c_{n-1} \xrightarrow{\alpha_n} c_n, x),
\end{multline}
where $x$ denotes an arbitrary element of the group $\mG(c_0)$.

In particular, the simplex $( c_0\xrightarrow{\alpha_1}c_1 \to \cdots \to c_{n-1} \xrightarrow{\alpha_n} c_n, x)$
will be degenerate if and only if at least one of its arrows is equal to $id$.
\end{proposition}

\subsection{Non-Abelian group diagram homology}

According to Proposition \ref{diagtosimp}, the simplicial group $C_*(\mC, \mG)$ can be used as a simplicial replacement of a group diagram. We note this in the following definition of non-abelian homology:

\begin{definition}\label{defhomol}
For an arbitrary small category $\mC$ and a group diagram $\mG: \mC\to {\Grp}$, the $n$th homologies of the category $\mC$ with coefficients in $\mG$ are the groups 
$\colim^{ \mC}_n \mG:=\pi_n(\coprod \mG)= \pi_n(C_*(\mC, \mG))$, for all $n\geq 0$.
\end{definition}

\section{Homotopy colimit for a non-abelian homology group}

In this section, we prove that the $n$th non-Abelian homology of a group diagram is isomorphic to the $(n+1)$th homotopy groups of the homotopy colimit for the diagram of the classifying spaces.

\subsection{Homotopy groups of a wedge for classifying spaces}

We will need a lemma on the permutability of the functor of a classifying space with coproducts up to homotopy.

Coproduct in the category of pointed simplicials
sets will be called a wedge. According to \cite[Ch.XII, \S 3.1]{bou1972}, the wedge is weakly equivalent to the homotopy colimit of a diagram over a discrete category in the category of pointed simplicial sets.

\begin{lemma}\label{bouquet}
For an arbitrary family of groups, there is a weak equivalence $\bigvee_{i\in I} BG_i \approx B (\MYstar_{i\in I}G_i)$.
\end{lemma}
{\sc Proof.}
This statement follows from a more general theorem of Whitehead mentioned in the papers
\cite[Theorem 4.0]{fie1984} and
\cite{bro1984}.
It will suffice for us to use the natural weak equivalence $B(G_1*G_2)\approx BG_1\vee BG_2$ following from the assertion proved in \cite[Lemma 3.8]{dwy1980}.

In the case of a finite set $I$ this assertion can be proved by induction.
If $I$ is infinite, then we use the weak equivalence $\bigvee_{i\in J} BG_i \approx B (\MYstar_{i\in J}G_i)$ for each finite subset of $J \subseteq I$.
Passing to a colimit with respect to a set of finite subsets and using the assertion \cite[Ch XII, \S 3.5]{bou1972} that the colimit of pointed simplicial sets with respect to a directed category is weakly equivalent to a homotopy colimit, we obtain the assertion of this lemma.
\hfill$\Box$


\subsection{Non-Abelian Homology and Homotopy Groups of a Homotopy Colimit}

We consider homotopy colimits in the category of pointed
simplicial sets. A homotopy colimit of classifying
spaces will be considered as a pointed simplicial set with a single vertex.

The following theorem
allows us to give a topological definition of the non-Abelian homology of a diagram
of a group.

\begin{theorem}\label{main2}
For an arbitrary group diagram $\mG: \mC\to {\Grp}$ and for all $n\geq 0$, there are natural isomorphisms
\begin{equation}\label{fmain2}
\colim^{\mC}_n \mG \cong
\pi_{n+1}(\hocolim^{\mC}B\mG).
\end{equation}
\end{theorem}
{\sc Proof.} By definition \ref{defhomol}, for all
$n\geq 0$, $\colim^{\mC}_n\mG= \pi_n(C_*(\mC, \mG))$.
Therefore, for all $n\geq 0$, we need to prove the existence of isomorphisms
\begin{equation}\label{fm3}
\pi_n(C_*(\mC, \mG)) \cong \pi_{n+1}(\hocolim^{\mC}B\mG)
\end{equation}

By Lemma \ref{spquill}, for every simplicial group
$G\in {\Grp}^{\Delta^{op}}$, there is a natural isomorphism
\begin{equation}\label{goja}
\pi_{n+1} \diag BG \cong \pi_n(G), \text{ for all $n\geq 0$},
\end{equation}
which we want to use.

Using Lemma \ref{bouquet}, for each $n\geq 0$, we arrive
at a weak equivalence

$$
(B\coprod \mG)_n = B\MYstar_{c_0\to \ldots \to c_n}\mG(c_0) \approx
\bigvee_{c_0\to \ldots \to c_n}B\mG(c_0)= (\coprod B\mG)_n.
$$
This weak equivalence will be natural in $[n]\in \Delta$. Consider pointed bisimplicial sets 
$B\coprod\mG$ and $\coprod B\mG$ as diagrams $\Delta^{op}\to \sSet_*$.
The homotopy colimit in $[n]\in \Delta^{op}$ leads to a weak equivalence 
$$
\hocolim^{\Delta^{op}} B\coprod\mG \approx \hocolim^{\Delta^{op}} \coprod B\mG.
$$
For every pointed bisimplicial set $Y\in \sSet^{\Delta^{op}}_*$ there exists a natural (in $Y$) weak equivalence $\hocolim^{\Delta^{op}}Y \approx \diag Y$ \cite[ch XII, \S 4.3]{bou1972}.
  This leads to a weak equivalence of the diagonals $\diag B\coprod\mG \approx \diag \coprod B\mG$.

Using this weak diagonal equivalence and substituting $G=\coprod\mG$ into the formula (\ref{goja}), we arrive at natural isomorphisms (with respect to $\mG$):
\begin{equation}\label{scheme}
\pi_{n+1} (\diag \coprod B\mG) \cong \pi_{n+1} (\diag B\coprod\mG) \cong \pi_n(\coprod \mG), \text{ for all $n\geq 0$}.
\end{equation}
By the simplicial replacement lemma \cite[Ch XII, 5.2]{bou1972}, for any diagram of pointed simplicial sets, the homotopy colimit $\hocolim^{\mC}Y$ is naturally weakly equivalent to $\diag \coprod Y$.
So the left side of the formula (\ref{scheme}) is equal to
$\pi_{n+1} (\hocolim^{\mC} B\mG)$.
By Definition \ref{defhomol} the right side of the isomorphism (\ref{scheme}) is equal to $\colim^{\mC}_n \mG$. Formula (\ref{fmain2}) received.
\hfill$\Box$

\begin{remark}
A similar formula in other terms is given by To\"en in \cite[Page 136]{toe2010}.
\end{remark}

\begin{remark}
To construct a pointed homotopy colimit using Thomason's theorem \cite[Theorem 1.2]{tho1979}, additional steps are needed, since Thomason's theorem has been proved for nonpointed homotopy colimits.
\end{remark}

\begin{corollary}\label{cor53}
For an arbitrary simplicial group $\mG\in {\Grp}^{\Delta^{op}}$ and
a non-negative integer $n$, we have an isomorphism
$$
\colim^{\Delta^{op}}_n \mG \cong \pi_n(\mG).
$$
\end{corollary}
{\sc Proof.} Consider an arbitrary integer $n\geq 0$.

In the case $\mC=\Delta^{op}$, Theorem \ref{main2} yields the following isomorphisms
$$
\colim^{\Delta^{op}}_n \mG\cong \pi_{n+1}\hocolim^{\Delta^{op}}B\mG.
$$
Now using the natural weak equivalence
$\hocolim^{\Delta^{op}}B\mG\to \diag B\mG$ \cite[Ch XII.3.4]{bou1972}
we get isomorphisms
$$
\pi_{n+1}\hocolim^{\Delta^{op}}B\mG\to \pi_{n+1}\diag B\mG
$$
for all $n+1\geq 0$.

Lemma \ref{spquill}, shows that $\pi_{n+1}(\diag B\mG) \cong \pi_n(\mG)$.

We obtain isomorphisms $\colim^{\Delta^{op}}_n \mG \cong \pi_n(\mG)$ for all $n\geq 0$.
\hfill$\Box$

For Abelian homology, this formula was proved by Bousfield and Kan \cite[XII.5.6]{bou1972},
in the case where $\mG$ is a simplicial Abelian group.
We have proved it for any simplicial group $\mG$.
\subsection{Abelian homology of diagrams over a free category}

Let $\cA: \mC\to \Ab$ be a diagram of abelian groups.
For the definition of $n$th Abelian homology $\coLim^{\mC}_n\cA$ one can take the homotopy groups $\pi_n(\coprod\cA)$ of the simplicial change of the diagram of Abelian groups $\cA$.
They will be isomorphic to the homology groups of the simplicial abelian group consisting of direct sums $C_n(\mC, \cA)= \oplus_{c_0\to \cdots \to c_n} \cA(c_0)$ with boundary operators $d^i_n$ and degeneracy operators $s^i_n$ acting on the elements of the terms of direct sums according to the formulas (\ref{fsimprep1})-(\ref{fsimprep2}) from the proposition \ref{diagtosimp}.
The Abelian homology groups are isomorphic to the homology groups of the complex of Abelian groups
$$
0 \xleftarrow{d_0} C_0(\mC, \cA) \xleftarrow{d_1} C_1(\mC, \cA) \xleftarrow{d_2}
\ldots \xleftarrow{d_n}C_n(\mC, \cA) \xleftarrow{d_{n+1}} \ldots
$$
whose differentials are $d_0=0$ and $d_n= \sum_{0\leq i\leq n}(-1)^n d^i_n$, for $n\geq 1$. Thus,
$$
\coLim^{\mC}_n\cA = \Ker d_n / \Imm d_{n+1}, \text{ for all } n\geq 0.
$$

Let $\mC= W\Gamma$ be the free category generated by the graph $\Gamma= (A\underset{s}{\overset{t}{\rightrightarrows}} V)$, where $A$ is the set of arrows (directed edges), $V$ is a set of vertices, $s: A\to V$ assigns each arrow its beginning, and $t: A\to V$ its end.

To calculate $\coLim^{W\Gamma}_1\cA$, \cite[Definition 2.1]{X2003} proposed the following

\begin{definition}
A flow on the graph $\Gamma= (A\underset{s}{\overset{t}{\rightrightarrows}} V)$ with coefficients in the diagram $\cA: W\Gamma \to \Ab$ is an element $f= ( f_{\gamma})_{\gamma \in A} \in \bigoplus\limits_{\gamma\in A}\cA(s(\gamma))$ such that for every vertex $v\in V$ there is the following ratio
$$
  \sum_{s(\gamma)\in v} f_{\gamma} = \sum_{t(\gamma)=v} \cA(\gamma)(f_{\gamma}).
$$
\end{definition}

The set of flows is a subgroup of the Abelian group $\bigoplus\limits_{\gamma\in A}\cA(s(\gamma))$.
It is proved \cite[Theorem 2.5]{X2003} that this flow subgroup is isomorphic to $\coLim^{W\Gamma}_1\cA$.
\begin{example}\label{pararrows}
Consider a graph consisting of two vertices $a$, $b$ and two arrows $a \xrightarrow{u_0} b$ and $a \xrightarrow{u_1} b$.
Denote by $\upuparrows$ the free category generated by this graph; it is obtained from the graph by adding the identical morphisms $1_a$ and $1_b$.
Let $\cA: \upuparrows \to \Ab$ be a digram of abelian groups.
The flow on the graph generating $\upuparrows$ is given by a pair of elements $f_{u_0}, f_{u_1}\in \cA(a)$ such that $f_{u_0}+f_{u_1}= 0$ and $\ cA(u_0)(f_{u_0})+\cA(u_1(f_{u_1})= 0$.
It follows that each flow is given by an arbitrary element $f_{u_0}\in \cA(a)$ such that $(\cA(u_0)- \cA(u_1))(f_{u_0})=0$. Hence $\coLim^{\upuparrows}_1 \cA= \Ker(\cA(u_0)-\cA(u_1))$.
\end{example}


\subsection{Abelian and non-Abelian homology of group diagrams}

Let us apply the non-Abelian and Abelian homology of group diagrams to develop a method for finding a non-zero homotopy group of least dimension for the homotopy colimit of classifying group spaces.

For the case when the colimit of a group diagram is equal to the trivial group, for diagrams defined on free
categories, we study the conditions under which the first non-Abelian and Abelian homology groups are isomorphic.

The content of this subsection is based on the following statement, which follows from 
Theorem  \ref{main2}.

\begin{corollary}\label{connect}
For an arbitrary small category $\mC$ and a diagram $\mG: \mC\to {\Grp}$ there is an isomorphism
\begin{equation}\label{farj}
\pi_1(\hocolim^{\mC} B\mG)\cong \colim^{\mC}\mG.
\end{equation}
In particular, a homotopy colimit of a diagram of classifying spaces is simply connected if and only if the colimit of a group diagram is equal to the trivial group.
\end{corollary}

\begin{example}
It would not be superfluous to check the formula (\ref{farj}) for a diagram that takes the values $\mG(c)= 0$, for all $c\in \Ob\mC$.
Since $B(0)=\star$, then according to Example \ref{forconst}, we get $hocolim^{\mC}(\Delta_{\mC}\star)= B\mC\times \star \diagup B\mC\times \star =\star$ and isomorphism (\ref{farj}).
\end{example}

Consider the problem of finding $n\geq 0$ such that $\colim^{\mC}_k\mG=0$ for all $0\leq k< n$ and $\colim^{\mC}_{n} \mG\not=0$.
  This number $n$ will be called the connectivity of the group diagram $\mG$ and denoted by $\cocon\mG$.
If $\colim_n^{\mC}\mG= 0$ for all $n\geq 0$, then we say that the diagram's connectivity is infinite and write $\cocon\mG= \infty$.
The connectivity of the diagram $\mG$ is equal to the largest $n$ for which the space $\hocolim^{\mC}B\mG$ is $n$-connected.

Let $\upuparrows$ be the category generated by the graph consisting of a pair of parallel arrows, considered in Example \ref{pararrows}.

\begin{example}\label{expar1}
Let $\mG: \upuparrows\to \Grp$ be a diagram taking the values $\mG(b)= S_3$ be the group of permutations of three elements, $\mG(a)= A_3$ be the subgroup of even permutations, $\mG( u_0): A_3\to S_3$ is an embedding of an alternating group into a symmetric group, $\mG(u_1): A_3 \to S_3$ is a homomorphism that takes constant values equal to the neutral element.

Since $\colim^{\upuparrows}\mG= \Coker\mG(u_0)= S_3/A_3= \ZZ_2$, we get that the group $\pi_1(\hocolim^{\mC} B\mG)= \ZZ_2$ is not trivial, and $\cocon\mG= 0$.

If in the considered diagram $\mG$ we replace the group $A_3$ by a cyclic subgroup of the second order $\ZZ_2\subset S_3$, then we get the equality 
$\colim^{\upuparrows}\mG= S_3/ \langle\ZZ_2 \rangle^ {S_3}=0$ and $\cocon\mG\geq 1$.
\end{example}

Before returning to  Example \ref{expar1}, we examine the conditions under which the values of the non-Abelian homology functor $\colim^{\mC}_1: \Grp^{\mC}\to \Grp$ are isomorphic
 to the values of the Abelian homology functor $\coLim ^{\mC}_1: \Ab^{\mC}\to \Ab$.

For this purpose we need the spectral sequence of Bousfield and Kan \cite[Ch. XII, \S5.7]{bou1972}, which exists for an arbitrary diagram of pointed sets $X(-): \mC\to \sSet$ and is located in the first quadrant
\begin{equation}\label{spbk}
E^2_{p,q}= \coLim^{\mC}_p \{\tilde{H}_q X(c)\}_{c\in\mC}
\Rightarrow \tilde{H}_{p+q}(\hocolim^{\mC} X),
\end{equation}
where $\tilde{H}_n$ denote the reduced homology groups of pointed simplicial sets \cite[Appendix II, \S 2]{gab1967}, $n\geq 0$.

\begin{proposition}\label{compar}
Let $\mG: W\Gamma\to \Grp$ be a group diagram over the free category $W\Gamma$ generated by an arbitrary graph in the sense of \cite[\S2.7]{mac2004}. We will assume that $\colim^{W\Gamma}\mG=0$.
  There is a natural (in $\mG\in \Grp^{W\Gamma}$) group homomorphism
  $$
  f: \colim^{W\Gamma}_1 \mG \to \coLim^{W\Gamma}_1\mG_{ab},
  $$
  where $\mG_{ab}$ is the Abelianization of the diagram $\mG$.
  This homomorphism $f$ is an isomorphism if and only if $\coLim^{W\Gamma}H_2(\mG)=0$.
  In particular, $f$ is an isomorphism if $\mG$ is a diagram of free groups.
\end{proposition}
{\sc Proof.} The non-Abelian homology group $\colim^{W\Gamma}_0\mG $ in dimension $n=0$ is isomorphic to $\colim^{W\Gamma}\mG$, and hence is equal to $0$ by assumption.

Let us calculate the first non-Abelian homology group of the diagram $\mG: W\Gamma\to \Grp$. Since $\colim^{W\Gamma}\mG= 0$, then $\pi_1(\hocolim^{W\Gamma} B\mG) = 0$ 
(Corollary   \ref{connect}).
It follows from the Hurewicz theorem that the group $\pi_2(\hocolim^{W\Gamma} B\mG)$ is isomorphic to the second integral homology group $H_2(\hocolim^{W\Gamma} B\mG)$.
To calculate this group, we use the spectral sequence (\ref{spbk}) for the diagram $X=B\mG$ of classifying spaces.
In this case, it will degenerate for all $n\geq 2$ into short exact sequences
\begin{equation}\label{shortexact}
0 \to \coLim^{W\Gamma} H_n (B\mG)
\to {H}_n (\hocolim^{W\Gamma} B\mG)
\to \coLim^{W\Gamma}_1 {H}_{n-1} (B\mG) \to 0
\end{equation}
where $\coLim^{W\Gamma}_1 {H}_{n-1} (B\mG)$ denotes the first Abelian homology group of the category $W\Gamma$ with coefficients in the diagram ${H}_{n-1 }(B\mG)$.

Since $H_2(\hocolim^{W\Gamma}B\mG)\cong \pi_2(\hocolim^{W\Gamma}B\mG) \cong colim^{W\Gamma}_1\mG$, and there are isomorphisms $H_1(B\mG(c))\cong \mG(c)_{ab}$,
  $H_2(B(\mG(c)))\cong H_2(\mG(c))$, then for $n=2$ the exact sequence (\ref{shortexact}) leads to the exact sequence natural in $\mG$
$$
0 \to \coLim^{W\Gamma} H_2 (\mG)
\to \colim^{W\Gamma}_1\mG
\to \coLim^{W\Gamma}_1 \mG_{ab} \to 0.
$$
It shows that the homomorphism $f: \colim^{W\Gamma}_1\mG
\to \coLim^{W\Gamma}_1 \mG_{ab}$ is an isomorphism if and only if $\coLim^{W\Gamma} H_2 (\mG)=0$.

In particular, if the diagram $\mG$ consists of free groups, then it can be given using the presentation $\langle E | R \rangle$, where $R=\{1\}$, whence, by Hopf's formula (see \cite[Theorem II.5.3]{broks1982}), the homology groups $H_2(\mG(c))$ are equal to $0$, which means that in this case the homomorphism $f$ is an isomorphism.
\hfill$\Box$

 The following example continues the group diagram exploration started in Example \ref{expar1}.

\begin{example}\label{expar2}
Our goal is to find the connectivity of the diagram
$$
\mG = ~
\left( \xymatrix{
\ZZ_2 \ar@<0.5ex>[r]^{<} \ar@<-0.5ex>[r]_{0} & S_3
}\right),
$$
which takes on objects the values $\mG(a)=\ZZ_2$, $\mG(b)= S_3$, and on morphisms $\mG(u_0): \ZZ_2\to S_3$ is the embedding of a cyclic subgroup of order $2$ into the group $S_3$, and $\mG(u_1)$ is a homomorphism onto a trivial subgroup.
We have established (Example \ref{expar1}) that $\colim^{\upuparrows}\mG=0$, whence the space $\hocolim^{\upuparrows}B\mG$ is simply connected. By the Hurewicz theorem, we get $\pi_2(\hocolim^{\upuparrows}B\mG)\cong
  H_2(\hocolim^{\upuparrows}B\mG)$.
 
Since $\upuparrows$ is a free category, then the sequence (\ref{shortexact}) will be exact, and for $n=2$ it will lead to the exact sequence

\begin{equation}\label{shortexact2}
0 \to \coLim^{\upuparrows} {H}_2 (\mG)
\to {H}_2 (\hocolim^{\upuparrows} B\mG)
\to \coLim^{\upuparrows}_1 {H}_{1} (\mG) \to 0
\end{equation}

Since $H_2(\ZZ_2)=0$ (see, for example, \cite[\S4.2]{kuz2006}) and $H_2(S_3)= 0$ \cite{mathoverflow180637}, then the diagram $H_2(\mG )$ is equal to the pair of trivial group homomorphisms $0 \rightrightarrows 0$, whence $\colim^{\upuparrows} {H}_2 (\mG)=0$.
The diagram $H_1(\mG)$ consists of two homomorphisms $\ZZ_2 \rightrightarrows (S_3)_{ab}=\ZZ_2$, one of which is the identity since $\mG(u_0)$ is a coretraction, and the other is equal to zero homomorphism $H_1(\mG(u_1))=\mG(u_1)_{ab}$,
because $\mG(u_1)$ is the null homomorphism.
Applying the method of calculating the values of the functor $\coLim^{\upuparrows}_1$ (example \ref{pararrows}), we obtain $\coLim^{\upuparrows}_1 H_1(\mG)=0$.
The exact sequence (\ref{shortexact2}) leads to $H_2(\hocolim^{\upuparrows} B\mG)= 0$.
  Hence $\colim^{\upuparrows}_1\mG= 0$, $\cocon\mG\geq 2$ (the space $\hocolim^{\upuparrows} B\mG)$ is $2$-connected) and
  $$
  \pi_3(\hocolim^{\upuparrows} B\mG))\cong H_3(\hocolim^{\upuparrows} B\mG))
  $$
  by the Hurewicz theorem.
  
 To calculate $H_3(\hocolim^{\upuparrows} B\mG)$ we have constructed an exact sequence (\ref{shortexact}), which for $n=3$ is equal to the exact sequence
  \begin{equation}\label{shortexact3}
0 \to \coLim^{\upuparrows} {H}_3 (\mG)
\to {H}_3 (\hocolim^{\upuparrows} B\mG)
\to \coLim^{\upuparrows}_1 {H}_{2} (\mG) \to 0
\end{equation}
We have shown that the diagram $H_2(\mG)$ consists of zeros, which implies that $\coLim^{\upuparrows}_1 {H}_{2} (\mG)=0$ and there is an isomorphism ${H} _3 (\hocolim^{\upuparrows} B\mG)\cong \coLim^{\upuparrows} {H}_3 (\mG)$. The diagram $H_3(\mG)$ consists of the Abelian groups $H_3(\mG(a))= H_3(\ZZ_2)\cong \ZZ_2$ \cite[\S4.2]{kuz2006} and $H_3(\mG( b))= H_3(S_3)\cong \ZZ_6$
  \cite{mathoverflow180637} between which two homomorphisms are given. Since the embedding $\mG(u_0): \ZZ_2 \to S_3$ is a coretraction, the homomorphism $H_3(\mG(u_0)): \ZZ_2\to \ZZ_6$ is a monomorphism.
  Since the zero homomorphism $\mG(u_1): \ZZ_2\to S_3$ is equal to the composition $\ZZ_2\to 0\to S_3$, then the homomorphism $H_3(\mG(u_1))$ passes through the group $H_3(0)= 0$ and hence $ H_3(\mG(u_1))$ will be the zero homomorphism. Hence the colimit $\coLim^{\upuparrows} {H}_3 (\mG)$ is isomorphic to the cokernel of the monomorphism $\ZZ_2 \to \ZZ_6$. Hence this colimit is isomorphic to the quotient  group $\ZZ_6/\ZZ_2 \cong \ZZ_3$.

We get the following answer: $\cocon\mG= 2$ and
$$
\colim^{\upuparrows}_2\mG= \pi_3(\hocolim^{\upuparrows} B\mG)= \ZZ_3.
$$
\end{example}

\begin{remark}
For the category $\mC$ having an initial object, the formula (\ref{farj}) can be obtained from \cite[Corollary 5.1]{far2004}.
\end{remark}
 
 \begin{example}\label{contr}
  Let's give an example to illustrate the formula (\ref{farj}).
  It shows that geometric intuition can be deceiving when constructing pointed homotopy 
  colimits.
  Let $\mC= \upuparrows$. Consider a diagram in the category
  $\sSet$ consisting of $X(a)= S^1$ is a simplicial circle (coequalizer of the pair 
  $\Delta(\partial^0), \Delta(\partial^1): \Delta[0]\rightrightarrows \Delta[1])$ \cite[II.2.5.2]{gab1967}, and let $X(b)=\Delta[0]$. Morphisms $X(u_0)=X(u_1): X(a)\to X(b)$ are defined as the only morphisms into a terminal object.

The geometric realization of the homotopy colimit of nonpointed simplicial sets $\hocolim^{\upuparrows}X$ consists of two cones over a circle whose vertices are identified.
Therefore, this homotopy colimit is homotopy equivalent to a sphere whose two points are identified. Its fundamental group is calculated using the Van Kampen theorem, as shown in \cite{brown_sphere2017}. It is equal to $\pi_1(\hocolim^{\upuparrows}X)= \ZZ$.
 
  Now consider the same diagram with only base vertices in $S^1$ and $\Delta[0]$.
 
  The construction of a homotopy colimit in $\sSet_*$ proceeds in steps \cite[XII.2]{bou1972}.
  The first step is to build a wedge of values on objects. In our case, this wedge will be equal to 
  $(S^1, \star)$ (the space $X(b)=\{\star\}$ will be identified with the base point).
  Then it is necessary to glue the reduced cylinders of mappings $X(u_0), X(u_1): X(a)\to X(b)$ to the resulting wedge (and in the general case the process continues further). In our case, these will be circles with a common border equal to the circle. The base points do not need to be connected, they are already identified. Therefore, the geometric realization of the pointed diagram is equal to the sphere. Whence the diagram of pointed simplicial sets satisfies  $\pi_1(\hocolim^{\upuparrows}X)= 0$.
 \end{example}
 
 \begin{remark} 
 The homotopy colimit of nonpointed simplicial sets considered in Example \ref{contr} is weak equivalent to $S^1\vee S^2$. To verify this, it suffices to replace the identified points by a segment.
 
  Farjoun proposed \cite{far2004} (see also \cite[Proposition 18.8.4]{hir2003}) a method for constructing the homotopy colimit of a diagram of pointed simplicial sets from non-pointed ones obtained by forgetting the property of points to be base. This homotopy colimit is equal to the quotient of the space of a nonpointed homotopy colimit with respect to the nerve of the category over which the diagram is given. Example \ref{contr} shows that this method is understandable.

  In our case, the classifying space of the category is a circle, and Farjoun's method yields $\hocolim^{\upuparrows}X\approx S^2$.
  \end{remark}

\section{Conclusion}

For an arbitrary small category $\mC$, we introduced non-Abelian homology of the category $\mC$ with coefficients in group diagrams and proved that they are isomorphic to derived functors of the functor $\colim^{\mC}: {\Grp}^{\mC}\to {\Grp} $.
We have proposed a formula relating non-Abelian homology to the homotopy groups of the homotopy colimit of classifying spaces. This formula provides a method for studying non-Abelian homology, and vice versa, there are applications of non-Abelian homology in homotopy theory.
We hope that this will contribute to the development of methods for calculating homotopy groups and generalization of the classical homology theory of small categories to non-Abelian homology
of group diagrams.

\end{document}